\documentclass[10pt,a4paper]{amsart}
\usepackage{amssymb,amsmath}

\usepackage{graphicx}

\newtheorem{theorem}{Theorem}[section]
\newtheorem{lemma}[theorem]{Lemma}

\theoremstyle{definition}

\newtheorem{remark}[theorem]{Remark}

\numberwithin{equation}{section}

\newcommand{\abs}[1]{\left|#1\right|}

\def\norm#1{\left\|#1\right\|}

\def \f12{\frac{1}{2}}

\newcommand{\eps}{\varepsilon}
\newcommand{\seq}[1]{\left\{#1\right\}}
\newcommand{\R}{\mathbb{R}}

\newcommand{\sgn}{\mathrm{sign}}
\newcommand{\Do}{\R\times\R_+}
\newcommand{\loc}{\mathrm{loc}}
\newcommand{\Test}{\mathcal{D}}
\newcommand{\Distr}{\mathcal{D}'}

\newcommand{\CM}{\mathcal{M}}

\newcommand{\pt}{\partial_t}

\newcommand{\ou}{\overline{u}}

\newcommand{\uenull}{u_{0}^{\delta}}
\newcommand{\ue}{u^{\eps,\delta}}
\newcommand{\we}{w^{\eps,\delta}}
\newcommand{\ve}{v^{\rho}}
\newcommand{\keps}{k^{\delta}}
\newcommand{\leps}{l^{\delta}}
\newcommand{\Feps}{F^{\rho}}
\newcommand{\Geps}{G^{\rho}}
\newcommand{\Heps}{H^{\rho}}
\newcommand{\Ieps}{I^{\eps,\delta}}

\newcommand{\Le}{\mathcal{L}^{\eps,\delta}}
\newcommand{\Leen}{\mathcal{L}_1^{\eps,\delta}}
\newcommand{\Leto}{\mathcal{L}_2^{\eps,\delta}}

\newcommand{\Letonull}{\mathcal{L}_{2,0}^{\eps,\delta}}
\newcommand{\Letoen}{\mathcal{L}_{2,1}^{\eps,\delta}}
\newcommand{\Letoto}{\mathcal{L}_{2,2}^{\eps,\delta}}
\newcommand{\Letotre}{\mathcal{L}_{2,3}^{\eps,\delta}}
\newcommand{\Letofire}{\mathcal{L}_{2,4}^{\eps,\delta}}
\newcommand{\Letofive}{\mathcal{L}_{2,5}^{\eps,\delta}}

\newcommand{\Hneg}{W^{-1,2}_{\mathrm{loc}}}

\newcommand{\Linf}{L^{\infty}}
\newcommand{\dx}{\, dx}
\newcommand{\dy}{\, dy}
\newcommand{\dt}{\, dt}
\newcommand{\dxi}{\, d\xi}

{%

\begin{enumerate}}%
{\end{enumerate}}

{%

\begin{enumerate}}%
{\end{enumerate}}

\begin{document}

\title[a resonant system of conservation laws]
{on the existence and compactness of\\ a two-dimensional 
resonant system of conservation laws}  

\author[K. H. Karlsen]{Kenneth H. Karlsen}
\address[Kenneth H. Karlsen]{\newline
         Centre of Mathematics for Applications \newline
         Department of Mathematics\newline
         University of Oslo\newline
         P.O. Box 1053, Blindern\newline
         N--0316 Oslo, Norway}
\email[]{kennethk@math.uib.no}
\urladdr{http://www.math.uio.no/\~{}kennethk/}
\author[M. Rascle]{Michel Rascle}
\address[Michel Rascle]{\newline Laboratoire 
J.A. Dieudonn\'e
\newline UMR CNRS n. 6621
\newline Universit\'e de Nice Sophia Antipolis,
\newline Parc Valrose, 06108 Nice Cedex 2, France.}
 \email[]{rascle@math.unice.fr}
\urladdr{http://math1.unice.fr/\~{}rascle/}
\author[E. Tadmor]{Eitan Tadmor}
\address[Eitan Tadmor]{\newline 
	Department of Mathematics and Institute for Physical Science and Technology\newline
	Center of Scientific Computation And Mathematical Modeling (CSCAMM)\newline
	University of Maryland\newline
	College Park, MD 20742 USA}
 \email[]{tadmor@cscamm.umd.edu}
\urladdr{http://www.cscamm.umd.edu/\~{}tadmor}

\date\today


\thanks{\textbf{Acknowledgment:} The research of K. H. Karlsen was supported by 
an Outstanding Young Investigators Award by the Research Council of Norway. 
The research of E. Tadmor was supported in part by NSF grant \#DMS04-07704 and ONR Grant 
\#N00014-91-J-1076} 

\begin{abstract}
We prove the existence of a weak solution to 
a two-dimensional resonant $3\times 3$ system 
of conservation laws with $BV$ initial data. 
Due to possible resonance (coinciding eigenvalues), spatial $BV$ estimates are 
in general not available. Instead, we use an entropy dissipation 
bound combined with the time translation invariance property 
of the system to prove existence based on a two-dimensional 
compensated compactness argument adapted from \cite{TRB}. 
Existence is proved under the assumption that the flux functions 
in the two directions are linearly independent. 
\end{abstract}

\maketitle

\section{Introduction} 
\label{Introduction}

This paper studies  certain two-dimensional resonant $3\times 3$ 
systems of conservation laws of the form
\begin{equation}
     \label{eq1}
     \begin{split}
          & k_t =0, \qquad l_t=0, \\
          & u_t + f(k,u)_x + g(l,u)_y=0,
     \end{split}
\end{equation}
which are augmented with $L^\infty\cap BV$ initial data 
\begin{equation}
     \label{data}
     k|_{t=0}=k(x,y), \quad 
     l|_{t=0}=l(x,y), \quad 
     u|_{t=0}=u_0(x,y).
\end{equation}
The goal is to prove that there exists a weak solution 
to \eqref{eq1}--\eqref{data}.

In recent years the one-dimensional 
version of the above system,
\begin{equation}
     \label{intro:eq1}
     \begin{split} 
          &k_t =0, \\ 
          &u_t + f(k,u)_x =0,
     \end{split}
\end{equation}
has received a considerable amount of attention. 
This system may be viewed as an alternative 
way  of writing a scalar conservation 
law with a discontinuous flux, namely
\begin{equation}
      \label{intro:eq2}
      u_{t} + f(k(x),u)_{x}=0.
\end{equation}

Equations like \eqref{intro:eq2}  occur in a variety of
applications, including flow in porous media,
sedimentation processes, traffic flow, radar shape-from-shading problems, 
blood flow,  and gas flow in a variable duct. 

If $k(x)$ is a smooth function, Kru\v{z}kov's
theory \cite{Kruzkov} tells us that there
exists a unique entropy solution to the initial
value problem for \eqref{intro:eq2}, for general  flux functions $f$. 
The scalar Kru\v{z}kov theory  does not apply when $k(x)$ is discontinuous. Instead it 
proves useful  to rewrite \eqref{intro:eq2} as a $2\times 2$ system 
of equations \eqref{intro:eq1}, which makes it possible 
to apply ideas from the theory of systems of 
conservation laws. 

As a starting point, it is necessary to introduce 
conditions on the flux $f(k,u)$ that guarantee
that solutions stay uniformly bounded. For example, one can require 
$f(k,a)=f(k,b)=0$ for all $k$, which in fact implies that 
the interval $[a,b]\subset \R$ becomes an invariant region.
The system \eqref{intro:eq2} has two eigenvalues, namely 
$\lambda_{1}=0$ and $\lambda_{2}=f_{u}(k,u)$. Consequently, if $f_{u}(k,u)$ vanishes 
for some value of $(k,u)$, then \eqref{intro:eq2} is  
nonstrictly hyperbolic and experiences  
so-called nonlinear resonant behavior, which 
implies that wave interactions are more 
complicated than in strictly hyperbolic 
systems. As a matter of fact, one cannot 
expect to bound the total variation of 
the conserved quantities directly, but only 
when measured under a certain singular mapping.  
A singular mapping that is relevant 
for \eqref{intro:eq1} is
$$
\Psi(k,u) = \int^u \abs{f_u(k,\xi)}\dxi.
$$
If $\seq{u^{\rho}}_{\rho>0}$ is a sequence of 
"reasonable" approximate solutions of \eqref{intro:eq1}, then 
one proves that the total variation of the transformed quantity 
$z^{\rho}:=\Psi(k,u^{\rho})$ 
is bounded independently of $\rho$.
Helly's theorem then gives convergence (along a subsequence) of
$z^{\rho}$ as $\rho \downarrow 0$. Since 
the continuous mapping $u \mapsto \Psi(k,u)$ is
one-to-one, $u^{\rho}$ also converges.

A singular mapping was used first by Temple 
\cite{Temple} to establish convergence of the
Glimm scheme (and thereby the existence of 
a weak solution) for a $2 \times 2$
resonant system of conservation laws modeling the displacement of oil
in a reservoir by water and polymer, which is now known to be
equivalent to a conservation law with a discontinuous
coefficient (see, e.g., \cite{KlingbroII}). Since then the 
singular mapping approach has been used and adapted by great many 
authors to prove existence of weak solutions  to resonant systems of 
conservation laws/scalar conservation laws with 
discontinuous flux functions,  by establishing 
convergence of various approximations 
schemes  (Glimm and Godunov schemes, front tracking, 
upwind and central type schemes, vanishing viscosity/smoothing method, $\ldots$), see (the 
list is far from being complete) 
\cite{Adimurthi:2005kx,Amadori:2004fk,Bachmann:2006fk,BBK, BKKR:C-T,BKRT-clarifier_II,
Gimsebro:Cauchy,JohnHong:2000,IsaacTemple:Resonance,IsTe:1992,IsaacTemple:Source,KlausenRisebro,KlingbroI,KlingbroII,
LTW:Comparison,LTW:Suppression,Noussair:SW,SegVov,TowersI,TowersII}.
Similar ideas have been used also in the context 
of degenerate parabolic equations \cite{KRT:FDM_degen}.  
Regarding uniqueness and entropy conditions for scalar 
conservation laws with discontinuous coefficients, see 
\cite{KRT:L1stable,Karlsen:2004sn} and the references therein.

As an alternative to the singular mapping approach, the 
papers \cite{KKR:Relax,KRT:CC,Karlsen:2004sn} 
has suggested to use the compensated compactness method 
and "scalar entropies" for the convergence analysis of approximate solutions. 
The results obtained with this 
approach are more general (and to some extent the proofs are easier) 
than those obtained with the singular mapping approach.

All the papers up to now have addressed the one-dimensional case. 
The aim of the present paper is to take a first look 
at the multi-dimensional case, which is completely unexplored.  
More precisely, we will prove the existence of at least one 
weak solution to the initial value problem for the two-dimensional 
system \eqref{eq1}. 

Our existence proof is based on studying the "$(\eps,\delta)\downarrow (0,0)$ limit" of 
classical solutions $\ue$ of the uniformly parabolic equation
\begin{align*}
      \ue_t + f(\keps,\ue)_x + g(\leps,\ue)_y=
      \eps\left(\ue_{xx}+ \ue_{yy}\right), 
      \qquad \eps>0,\, \delta>0,
\end {align*}
where $\keps, \leps$ converge to $k,l$ in 
$L^1_{\loc}(\R^2)$, respectively, as $\delta\downarrow 0$.  

Observe that  we are essentially considering a scalar approximation scheme for \eqref{eq1}, see 
\cite{BKKR:C-T,BKRT-clarifier_II,KRT:FDM_degen,KKR:Relax,KRT:CC,Karlsen:2004sn,TowersI,TowersII} for other scalar approximation schemes for  one-dimensional discontinuous flux problems.

Although spatial $BV$ bounds are out of reach, 
we still have a time translation invariance property at our 
disposal, which, together with the assumption of $BV$ initial data, implies 
that $\ue_{t}$ is uniformly bounded in $L^1$. 
Consider three functions $F(k,u)$, $G(l,u)$, $H(k,l,u)$ defined by 
$$
F_{u}=(f_{u})^2, \quad G_{u}=(g_{u})^2,\quad 
H_{u}=f_{u}g_{u}.
$$
We prove, at least under the assumption that 
$\eps$ and $\delta$ are of 
comparable size, that the two sequences
$$
F(k(x,y),\ue)_{x} +H(k(x,y),l(x,y),\ue)_{y}  
$$
and 
$$
H(k(x,y),l(x,y),\ue)_{x} + G(k(x,y),\ue)_{y}   
$$
are compact in $\Hneg(\R^2)$, for each fixed $t>0$.  

The crux of the convergence analysis is then to  prove 
that the above  $\Hneg(\R^2)$ compactness is sufficient  
to establish a "two-dimensional" compensated compactness argument 
in the spirit of the classical Tartar-Murat results for one-dimensional 
conservation laws \cite{MuratII,MuratI,Murat:Hneg,TartarI,TartarII} 
(see also \cite{Benzoni-Gavage:1994fk}). 
Here we follow the recent two-dimensional compensated compactness framework 
developed in Tadmor et. al. \cite{TRB} for \emph{nonlinear} 
conservation laws. We extend their results to the case involving 
additional discontinuous "variable coefficients". Accordingly, we make 
the \emph{nonlinearity assumption} that 
for each fixed $k,l$ the functions $u\mapsto f_u(k,u)$ and  $u\mapsto g_{u}(l,u)$ are 
almost everywhere linearly independent (see (\ref{ass:gen_nonlin}) in 
the next section for a precise statement).
Our main existence result is based on an application of the two-dimensional 
compensated compactness  lemma with  "variable coefficients" --- lemma \ref{lem:CC} stated in 
Section \ref{sec:CC} below. Granted the nonlinearity assumption, it then 
yields that (a  subsequence of) $\ue(\cdot,\cdot,t)$ converges 
in $L^1_{\loc}(\R^2)$ to a bounded function $u(\cdot,\cdot,t)$, for a.e.~$t>0$. 
Since $\ue$ is uniformly $L^1$ Lipschitz continuous in time we 
obtain, in Section \ref{sec:proof} below, our main  Theorem \ref{thm:main}, stating
that $\ue\to u$ in $L^1_{\loc}(\R^2\times\R_{+})$ 
and that the limit function $u$ is a weak solution of \eqref{eq1}--\eqref{data}.

Although we have chosen to analyze the 
vanishing viscosity/smoothing method, the techniques 
used here for that purpose can also be applied to 
various numerical schemes, including appropriate two-dimensional versions of the 
scalar finite difference schemes studied in  \cite{KKR:Relax,Karlsen:2004sn,TowersI,TowersII}. 

 \section{Assumptions and statement of main results}
\label{sec:results}
We start by listing the assumptions on the initial conditions $u_0$ and the fluxes $k,l,f,g$ 
that are needed for the existence result.

Regarding the initial function we assume
\begin{equation}
   \label{ass:init_data_en}
   u_0 \in \Linf(\R^2)\cap BV(\R^2), \qquad a \leq u_0
   \leq b \quad \text{for a.e.~in $\R^2$}.
\end{equation}

For the discontinuous coefficients $k,l:\R^2\to \R$ we assume 
\begin{equation}
   \label{ass:kl}
   \begin{cases}
         k,l \in \Linf(\R^2)\cap BV(\R^2), \\
         \alpha \leq k,l
        \leq \beta \quad \text{a.e.~in $\R^2$}.
   \end{cases} 
\end{equation}

For the flux functions $f,g:[\alpha,\beta]\times [a,b]\to \R$ we assume 
\begin{equation}
   \label{ass:fg_en}
   \begin{cases}
      \text{$u\mapsto f(k,u),u\mapsto g(l,u)\in C^2[a,b]$ for all $k,l\in [\alpha,\beta]$};\\
      \text{$k\mapsto f(k,u),l\mapsto g(l,u)\in C^1[\alpha,\beta]$ for all $u\in [a,b]$}.
   \end{cases}
\end{equation}

Moreover, we make the \emph{nonlinearity assumption} which excludes the possibility of 
$\xi_1 f(k,u)+\xi_2 g(l,u)$ being an affine function (in $u$) on any nontrivial interval for all 
$k,l\in [\alpha,\beta]$,
\begin{subequations}\label{ass:gen_nonlin}
\begin{equation*}
\forall  \abs{\xi}=1  
         \text{ and }~k,l\in [\alpha,\beta]: \quad \xi_1 f(k,\cdot)+\xi_2 g(l,\cdot) 
         \equiv\!\!\!\!\!/  \ \ \text{affine function on  any nontrivial interval}.	
\end{equation*}
In its slightly stronger version, this assumption requires 
that $f_u(k,\cdot)$ and $g_u(l,\cdot)$ are a.e. linearly independent so that the symbol 
$s(\xi,k,l, u):=\xi_1 f_{u}(k,u)+\xi_2 g_{u}(l,u)$ satisfies 
\begin{equation*}
\forall  \abs{\xi}=1
    \     \text{ and }~k,l\in [\alpha,\beta]: \qquad 
	{\rm meas}\{ u \ | \ s(\xi,k,l,u)=0\} = 0.	
\end{equation*}
\end{subequations}
This is a straightforward generalization of the notion of nonlinearity found in \cite{LPT}, in their 
study of kinetic formulations for nonlinear conservation laws.

Finally, we need to know that our approximate solutions
stay uniformly bounded. For example, this is ensured by 
the assumption
\begin{equation}
   \label{ass:fg_to}
   \text{$f(k,a),f(k,b),g(l,a),g(l,b)=0$ for all $k,l\in [\alpha,\beta]$},
\end{equation}
which implies that the interval $[a,b]$ becomes 
an invariant region. Of course, one can relax assumption \eqref{ass:fg_to}. A sufficient condition for 
the invariance of the interval $[a, b]$ is that the divergence of the 
vector field $(x, y) \mapsto (f(k(x, y),u), g(l(x, y), u))$ is nonnegative 
when $u = b$ and nonpositive when $u = a$. Let us emphasize 
that an assumption like \eqref{ass:fg_to} is essential to our analysis; Without it 
solutions can possess concentration effects, which is a well-known 
feature of, for example, linear transport equations with discontinuous coefficients.

We are now ready to state our main result.
\begin{theorem}
\label{thm:main}
Suppose \eqref{ass:init_data_en}, \eqref{ass:kl}, \eqref{ass:fg_en}, \eqref{ass:gen_nonlin}, 
and  \eqref{ass:fg_to} hold. Then, there exists 
a weak solution of the initial value problem \eqref{eq1}--\eqref{data},
$u\in  \Linf(\R^2\times\R_{+})\cap \mathrm{Lip}(\R_{+};L^1(\R^2))$,  
 satisfying 
\begin{equation*}
      \begin{split}
            &\int_{\R_+}\int_\R \Bigl(u\phi_t + f(k(x,y),u)\phi_x
            + g(l(x,y),u)\phi_y\Bigr)\dx\dt\\
            &\qquad            + \int_\R u_0(x)\phi(x,0)\dx=0, 
            \qquad \forall \phi\in \Test(\R^2\times[0,\infty)).
      \end{split}
\end{equation*}
The weak solution, $u$, can be constructed as 
a strong $L^1_{\loc}(\R^2\times\R_{+})$-limit of classical solutions 
$\ue$ of  uniformly parabolic problems,  
\begin{equation}
	\label{eq:L1contr-eqn}
	\ue_t + f(\keps,\ue)_x + g(\leps,\ue)_y=   \eps\Delta \ue,
\end{equation}
with the smoothly mollified coefficients, $\keps:=\omega_\delta\star k$ 
and $\leps:=\omega_\delta\star l$ (outlined in section \ref{thm:main} below).
\end{theorem}

The proof of this theorem is given 
in the following two sections. Remark that the $BV$ assumption on the coefficients $k,l$ made in \eqref{ass:kl} 
is used twice in this paper. First, it is used to prove Lipschitz regularity in time, 
in lemma \ref{lem:time} below; then,  we use it to prove 
$\Hneg(\R^2)$ compactness of the entropy production for 
each fixed $t>0$ in lemma \ref{lem:Hneg_viscous} below.

\noindent
We close this section with the following summary. 

\begin{remark}\normalfont
Stated differently, Theorem \ref{thm:main} shows
that there exists a weak solution to the 
following two-dimensional scalar conservation law 
with discontinuous coefficients 
$k,l\in L^\infty(\R^2)\cap BV(\R^2)$:
\begin{align*}
      &u_t+f(k(x,y),u)_x + g(l(x,y),u)_y=0, \\
      & u|_{t=0}=u_{0}\in  L^\infty(\R^2)\cap BV(\R^2).
\end{align*}
We note in passing that the solution operator in this case of 
discontinuous "variable coefficients"  is not translation invariant in space 
and hence the $L^1$-contraction property of \eqref{eq:L1contr-eqn} does not 
imply spatial $BV$ compactness. 

Moreover, if we let $u_0(\cdot) \mapsto u(t,\cdot)$ denote the 
mapping of \eqref{eq1}--\eqref{data}, so that $u(t,\cdot)$ is a
(vanishing viscosity) weak solution constructed in theorem \ref{thm:main}, then by 
adapting standard arguments, we can prove that the mapping is compact 
with respect to the $L^1_\loc$ norm. 

\end{remark}

\section{A compensated compactness lemma}
\label{sec:CC}
In this section we prove a "two-dimensional" compensated compactness
lemma. We refer  \cite{Chen:LN,Lu:Book,MuratII,MuratI,Murat:Hneg,TartarI,TartarII} 
for background information on the compensated compactness theory.
We start by recalling the celebrated div-curl lemma.

\begin{lemma}[div-curl]
\label{lem:divcurl}
Let $\Omega\subset \R^2$ be an open domain and let 
$\rho>0$ denote a parameter taking its values in a 
sequence which tends to zero. Suppose $D^\rho \rightharpoonup  \overline{D}$, 
$E^\rho \rightharpoonup \overline{E}$ in $\bigl(L^2(\Omega)\bigr)^2$ 
and $ \left\{\mathrm{div}\, D^\rho\right\}_{\rho>0}$, $\left\{\mathrm{curl}\, E^\rho\right\}_{\rho>0}$ 
belong to a compact subset of $\Hneg(\Omega)$. Then, after extracting 
a subsequence if necessary, we have 
$D^\rho\cdot E^\rho \to \overline{D}\cdot \overline{E}$ 
in $\Distr(\Omega)$ as $\rho \downarrow 0$.
\end{lemma}

The compensated compactness lemma below 
is tailored for two-dimensional equations, 
whose spatial part involve discontinuous coefficients:
$$
f(k(x,y),v(x,y))_{x} + g(l(x,y),v(x,y))_{y}.
$$
If $g(l,u)=g(u)$ and $f(k,u)=f(u)$, 
then the lemma below coincides with the two-dimensional 
result of \cite[Theorem 3.1]{TRB}. If we set $g=0$ then 
the result coincides with Tartar's compensated 
compactness lemma for the one-dimensional scalar
conservation with genuinely nonlinear flux $f$.

\begin{lemma}[Compensated compactness]
\label{lem:CC}
Let $\Omega\subset \R^2$ be an open domain. 
Let $k,l,f,g$ be functions satisfying \eqref{ass:kl}, \eqref{ass:fg_en}, \eqref{ass:gen_nonlin}, 
and  \eqref{ass:fg_to}.  Suppose $\seq{\ve(x,y)}_{\rho>0}$ is
a sequence of measurable functions on $\Omega$ that satisfies
the following two conditions:
\begin{enumerate}
   \item There exist two finite constants $a<b$
    independent of $\rho$ such that
    $$
    a\le \ve(x,y)\le b \quad \text{for a.e.~$(x,y)\in \Omega$}.
    $$

   \item Let the functions $F,G,H$ be defined by 
   $$
   F_u(k,u)=\left(f_u(k,u)\right)^2,\quad 
   G_u(l,u)=\left(g_u(l,u)\right)^2, \quad 
   H_u(k,l,u)= f_u(k,u)g_u(l,u).
   $$
    We assume that  the two sequences 
    \begin{equation}
       \label{assumption:Hneg:2D}
       \begin{split}
          &\Bigl\{F\left(k(x,y),\ve\right)_x
          + H\left(k(x,y),l(x,y),\ve\right)_y\Bigr\}_{\rho>0},\\ 
          &\Bigl\{H\left(k(x,y),l(x,y),\ve\right)_x
          + G\left(l(x,y),\ve\right)_y\Bigr\}_{\rho>0}
       \end{split}
    \end{equation}
    belong to a compact subset of $\Hneg(\Do)$. 
\end{enumerate}
Then, there exists a subsequence of $\seq{\ve(x,y)}_{\rho>0}$
that converges a.e.~to a function $v\in \Linf(\R^2)$, and $a\le v(x,y)\le b$ 
for a.e.~$(x,y)\in \Omega$.
\end{lemma}

\begin{proof}
To simplify the notation let
\begin{align*}
      &\Feps :=F\left(k(x,y),\ve\right), \quad \Geps :=G\left(l(x,y),\ve\right), \quad
      \Heps :=H\left(k(x,y),l(x,y),\ve\right),
\end{align*}
and denote their $L^\infty(\Do)$ weak-$\star$ limits 
by $\overline{F},\overline{G},\overline{H}$, respectively. 
Introduce the vector fields
$$
D^\rho = \Bigl(\Feps,\Heps\Bigr),
\qquad 
E^\rho = \Bigl(-\Geps,\Heps\Bigr),
$$
and denote their respective $L^\infty(\Do)$ weak-$\star$ limits by $\overline{D},\overline{E}$.

Thanks to \eqref{assumption:Hneg:2D}, we can apply the div-curl lemma to the 
sequences $\seq{D^\rho}_{\rho>0}$, $\seq{E^\rho}_{\rho>0}$ to produce
$$
\overline{D\cdot E}=
\overline{D}\cdot \overline{E} 
\quad \text{a.e.~in $\Omega$},
$$
that is,
$$
\overline{H^2-FG}
=\left(\overline{H}\right)^2 -\overline{F}\,\overline{G},
$$
which implies
\begin{equation}
   \label{eq:MTrelation:2D_I}
   \begin{split}
      &\overline{\left(H-\overline{H}\right)^2
      - \left(F -\overline{F}\right)\left(G-\overline{G}\right)}
      =0.
   \end{split}
\end{equation}

Fix $c=c(x,y)\in L^\infty(\Omega)$. Following \cite{TRB}, we now consider the function 
$I:[a,b]\to \R$ defined by 
\begin{align*}
   I(v)&=\Big(H\left(k(x,y),l(x,y),v\right)-H\left(k(x,y),l(x,y),c\right)\Big)^2
   \\ & \quad 
    - \Big(F\left(k(x,y),v\right) -F\left(k(x,y),c\right)\Big)\cdot
   \Big(G\left(l(x,y),v\right)-G\left(l(x,y),c\right)\Big).
\end{align*}

Note that 
\begin{align*}
   I(\ve) &=
   \Bigl(\left[\,\Heps-\overline{H}\,\right] + 
   \left[\,\overline{H}-H\left(k(x,y),l(x,y),c\right)\,\right]\Bigr)^2
   \\ & \qquad 
    - \Bigl(\left[\,\Feps -\overline{F}\,\right] 
    + \left[\overline{F}-F\left(k(x,y),c\right)\right]\Bigr)\cdot
       \Bigl(\left[\,\Geps- \overline{G}\,\right]
    +\left[\,\overline{G}-G\left(l(x,y),c\right)\,\right]\Bigr).
\end{align*}

Using this and \eqref{eq:MTrelation:2D_I}, we compute 
\begin{equation}
   \label{eq:MTrelation:2D_II}
        \overline{I(v)}
          = \Bigl(\overline{H}-H\left(k(x,y),l(x,y),c\right)\Bigr)^2
                    - \Bigl(\overline{F} -F\left(k(x,y),c\right)\Bigl)\cdot
         \Bigl(\overline{G}-G\left(l(x,y),c\right)\Bigr).
\end{equation}

By the Cauchy-Schwarz inequality we have 
for any $u\in [a,b]$
\begin{equation}
   \label{eq:CS:2D}
   \begin{split}
      &\Bigl(H\left(k(x,y),l(x,y),u\right)-H\left(k(x,y),l(x,y),c\right)\Bigr)^2 
      \\
      & \qquad =\left(\int_{c}^{u} f_{u}(k(x,y),\xi)g_{u}(l(x,y),\xi)\dxi\right)^2
      \\
       & \qquad \le \Bigl(F\left(k(x,y),u\right)-F\left(k(x,y),c\right)\Bigl)
      \Bigl(G\left(l(x,y),u\right)-G\left(l(x,y),c\right)\Bigr),
   \end{split}
\end{equation}
and hence $I(\cdot)\le 0$ with $I(c)=0$. Thanks to  \eqref{ass:gen_nonlin}, the 
Cauchy-Schwarz inequality in (\ref{eq:CS:2D})  is in fact 
a \emph{strict} inequality. This shows that the function 
 $I(v)$ has a strict global maximum at $v=c$. 

Since $u\mapsto F\left(k(x,y),u\right)$ is strictly 
increasing, we can choose $c$ as
$$
c(x,y):=F^{-1}\left(k(x,y),\overline{F}(x,y)\right),
$$
so that \eqref{eq:MTrelation:2D_II} becomes
\begin{equation*}
   \overline{I(v)}
   = \Bigl(\overline{H}-H\left(k(x,y),l(x,y),c\right)\Bigr)^2.
\end{equation*}

Since $I(\cdot)\le 0$, we conclude that 
$\overline{H}=H\left(k(x,y),l(x,y),c\right)$, and thus 
$\overline{I(v)}=0 $. In fact, we have $I(\ve)\to 0$ a.e.~in $\Omega$.

Using the fact that $c$ is a strict 
maximizer of $I(v)$, we have
$$
I(v)\le -C_\alpha,
\quad \text{whenever $\abs{v- c}>\alpha$},
$$
for some constant $C_\alpha>0$ that depends on $\alpha$.
Consequently,
\begin{equation*}
      \mathrm{meas}\left\{\abs{\ve - \ou}>\alpha\right\}
      \le \frac{-1}{C_\alpha} \iint_{\Omega\cap \abs{\ve - c}>\alpha}
      I(\ve(x,y))\dx\dy\to 0 \quad \text{as $\rho\downarrow 0$}.
\end{equation*}
Since $\alpha>0$ was arbitrary, $\ve\to c$ in measure, which in turn implies 
that a subsequence of $\seq{\ve}_{\delta>0}$ converges to $c$ a.e.~in $\Omega$.
\end{proof}

We remark that the idea of using the Cauchy-Schwarz inequality along 
the lines of \eqref{eq:CS:2D} for proving strong compactness 
can be traced back to \cite{Rascle:vn,Rascle:1984kx}

\section{Proof of Theorem \ref{thm:main}} 
\label{sec:proof}

Let $\keps,\leps, \uenull$ be smooth functions converging 
strongly to $k,l, u_0$ respectively. More precisely, let 
$\omega_{\delta} \in C_0^{\infty}(\R)$ 
be a nonnegative function satisfying
$$
\omega(x)\equiv 0\quad
\text{for $|x|\geq 1$}, \qquad 
\iint_{\R^2}\omega(x)\,dx=1.
$$
For $\delta>0$, let $\omega_{\delta}(x)=
\frac{1}{\delta^2}\omega\left(\frac{x}{\delta}\right)$ and
introduce the mollified functions
$$
\keps=\omega_{\delta}\star k, \quad 
\leps=\omega_{\delta}\star l.
$$
We approximate the initial data $u_0$ by 
cut-off and mollification as follows:
$$
\uenull =\omega_{\delta} \star (u_0\chi_{\mu}),
$$
where $\chi_{\delta}(x)=1$ for $\abs{x}\leq 1/\delta$ and 0 otherwise. In
particular, we have the estimate 
\begin{equation}
	\label{eq:approx-property}
	\norm{\left(\uenull\right)_{xx}+\left(\uenull\right)_{yy}}_{L^1(\R^2)}
	\leq \frac{1}{\delta}\iint_{\R^2}
	\left( \abs{\left(\uenull\right)_x}
	+\abs{\left(\uenull\right)_y}\right)\dx\dy
	\leq\frac{1}{\delta}\abs{u_0}_{BV(\R^2)}.
\end{equation}
Observe that for $h^\delta=\keps,\leps,\uenull$ and 
$h=k,l,u_{0}$, we have $h^{\delta}\in C^\infty(\R^2)$ and 
$$
h^\delta\to h \quad \text{a.e.~in $\R^2$ and in $L^p(\R^2)$ 
for any $p\in [1,\infty)$ as $\delta\downarrow 0$}.
$$
Additionally, $\uenull$ is compactly supported.

Now let $\ue$ be the unique classical solution of the uniformly 
parabolic equation
\begin{equation}
      \label{eq:viscous}
      \ue_t + f(\keps,\ue)_x + g(\leps,\ue)_y=
      \eps\left(\ue_{xx}+ \ue_{yy}\right),
\end{equation}
with initial data $\ue|_{t=0}=\uenull$. The proof proceeds through a series of lemmas, 
which in the end show that for each $t\in [0,T]$ a subsequence 
of $\ue(\cdot,\cdot,t)$ converges a.e.~as $\varepsilon,\delta \downarrow 0$.

Our first lemma confirms the uniform bound. 
\begin{lemma}[$L^\infty$ bound]\label{lem:max}
There exists a constant $C>0$, independent of $\eps,\delta$, such that
$$
\norm{\ue(\cdot,\cdot,t)}_{L^{\infty}(\R^2)}\le C,
\qquad \text{for all}\quad t\in (0,T).
 $$
\end{lemma}

\begin{proof}
The proof is standard and 
exploits assumption \eqref{ass:fg_to} to  conclude that 
$a\le \ue(x,y,t)\le b$ for a.e.~$(x,y,t)\in \R^2\times\R_{+}$.
\end{proof}

Using that \eqref{eq:viscous} is 
translation invariant in time, we can 
prove that $\ue_{t}$ is uniformly 
bounded in $L^\infty(\R_{+};L^1(\R^2))$.

\begin{lemma}[Lipschitz regularity in time]\label{lem:time}
Suppose the two smoothing parameters $\eps$ and $\delta$ are kept 
in balance in the sense that 
\begin{equation}
   \label{eq:balance}
   \delta = C \eps, \quad \text{for some constant $C>0$.} 
\end{equation}
There is a constant $C_0$ (which is possibly dependent on $u_0$ but otherwise 
is independent of $\eps,\delta$), such that for any $t>0$
$$
\iint_{\R^2}\abs{\pt \ue (\cdot,\cdot,t)}\dx\dy\le C_0.
$$
\end{lemma}

\begin{proof}
To prove this, set $\we=\ue_t$. Then $w$ satisfies
$$
\we_t + (f_u(\keps,\ue)\we)_x + (g_u(\leps,\ue)\we)_x=
\eps\left(\we_{xx} + \we_{yy}\right).
$$
Multiplying by $\sgn(\we)$ gives, in the sense of distributions,
\begin{align*}
      &\abs{\we}_t + (f_u(\keps,\ue)\abs{\we})_x + (g_u(\leps,\ue)\abs{\we})_y
      \\
      & \qquad 
      =\eps\left(\abs{\we}_{xx} + \abs{\we}_{yy}\right)
      - \eps\sgn'(\we)\left(\left(\we_x\right)^2+\left(\we_y\right)^2\right),
\end{align*}
since $f_u(\keps,\ue)\we \sgn'(\we)\we_x = 
g_u(\keps,\ue)\we \sgn'(\we)\we_y =0$.
Hence
$$
\frac{d}{dt} \iint_{\R^2} \abs{\we}(x,y,\cdot)\dx\dy\le 0,
$$
which, due to \eqref{eq:approx-property} and \eqref {eq:balance}, concludes 
the proof.
\end{proof}

Thanks to the previous lemma, we also have 
uniform $L^{\infty}(\R_{+};L^2(\R^2))$ control 
over $\sqrt{\eps}\ue_{x}$ and  $\sqrt{\eps}\ue_{y}$. 

\begin{lemma}[Entropy dissipation bound]
\label{lem:energy}
There is a constant $C$, independent of $\eps,\delta$, such that 
$$
\eps\iint_{\R^2} \Bigl( \left(\ue_x (\cdot,\cdot,t)\right)^2
+\left(\ue_y (\cdot,\cdot,t)\right)^2\Bigr)\dx\dy\le C, 
\qquad \text{for any $t>0$.}
$$
\end{lemma}

\begin{proof}
Multiplying \eqref{eq:viscous} by 
$\ue$ and then integrating yield
\begin{align*}
   &\iint_{\R^2}\eps \Bigl( \left(\ue_x(\cdot,\cdot,t)\right)^2
   +\left(\ue_y(\cdot,\cdot,t)\right)^2\Bigr)\dx\dy
   \\ & \quad = - \iint_{\R^2} u \ue_t \dx\dy
   \\ & \qquad 
   + \iint_{\R^2}\left[ \left(\int_0^{\ue} f_u(\keps,\xi)\dxi\right)_x
   - \left(\int_0^{\ue} f_{k}(\keps,\xi)\dxi\right)
   \keps_x(x,y)\right]\dx\dy
   \\ & \qquad 
   + \iint_{\R^2}\left[ \left(\int_0^{\ue} g_u(\leps,\xi)\dxi\right)_y 
   - \left(\int_0^{\ue} g_{l}(\leps,\xi)\dxi\right)\leps_y\right]\dx\dy.
\end{align*}

In view of Lemmas \ref{lem:max} and \ref{lem:time} and 
the $BV$ regularity of the coefficients, we derive easily 
the uniform bound
\begin{align*}
     &\eps\iint_{\R^2} \Bigl( \left(\ue_x (\cdot,\cdot,t)\right)^2
     +\left(\ue_y (\cdot,\cdot,t)\right)^2\Bigr)\dx\dy
     \\
     & \le C\left(\norm{\ue_{t}}_{L^\infty(\R_{+};L^1(\R^2))} 
     +\abs{k}_{BV(\R^2)} +\abs{l}_{BV(\R^2) }\right),\qquad t>0,
\end{align*}
for some constant $C$ that is dependent on $\|\ue\|_{L^\infty(\R^2)}$ but 
otherwise is independent of $\eps,\delta$.
\end{proof}

\begin{lemma}[Pre-compactness at each time instant]
\label{lem:Hneg_viscous}
Suppose the two parameters $\eps$ and $\delta$ are kept 
in balance in the sense that \eqref{eq:balance} holds. 
With $F$, $G$, and $H$ defined in Lemma \ref{lem:CC}, 
the two sequences 
\begin{align*}
   &\Bigl\{F\left(k(x,y),\ue\right)_x
   + H\left(k(x,y),l(x,y),\ue\right)_y\Bigr\}_{\eps,\delta>0},\\ 
   &\Bigl\{H\left(k(x,y),l(x,y),\ue\right)_x
   + G\left(l(x,y),\ue\right)_y\Bigr\}_{\eps,\delta>0}
\end{align*}
then belong to a compact subset of $\Hneg(\R^2)$, for each fixed $t>0$.
\end{lemma}

\begin{proof}
Let $\phi=\phi(x,y)\in \Test(\R^2)$, and, for each fixed $t>0$, introduce 
the distribution 
$$
\left\langle \Le, \phi \right\rangle
=\iint_{\R^2} \Bigl(F\left(k(x,y),\ue\right)\phi_x
+ H\left(k(x,y),l(x,y)\ue\right) \phi_y\Bigr)\dx\dy.
$$
Let us first write $\Le = \Leen+\Leto$, where
\begin{align*}
   \left\langle \Leen, \phi \right\rangle
   &=\iint_{\R^2} \Bigl(F\left(k(x,y),\ue\right)
   -F\left(\keps(x,y),\ue\right)\Bigr)\phi_x\dx\dy
   \\ & \qquad
   +\iint_{\R^2}\Bigl(H\left(k(x,y),l(x,y),\ue\right)
   -H\left(\keps(x,y),\leps(x,y),\ue\right)\Bigr)\phi_y\dx\dy,\\
   \left\langle \Leto, \phi \right\rangle
   &=\iint_{\R^2} \Bigl(F\left(\keps(x,y),\ue\right)\phi_x
   + H\left(\keps(x,y),\leps(x,y), \ue\right) \phi_y\Bigr)\dx\dy.
\end{align*}

In what follows, we let $\Omega$ denote an
arbitrary but fixed bounded open subset of $\R^2$. 
With $\phi \in W^{1,2}_0(\Omega)$, we 
have by H\"older's inequality
$$
\left|\left\langle \Leen, \phi \right\rangle\right|
\le C \left(\norm{k - \keps}_{L^2(\Omega)} +
\norm{k - \keps}_{L^2(\Omega)}\right)
\norm{\phi}_{W^{1,2}_0(\Omega)}
\to 0,
$$
as $\delta\downarrow 0$. Thus, $\left\{\Leen\right\}_{\eps,\delta>0}$ 
is compact in $W^{-1,2}(\Omega)$, for each fixed $t$.

Multiplying \eqref{eq:viscous} by 
$f_u(\keps(x,y),\ue)$ yields
\begin{align*}
   &f(\keps,\ue)_t 
   + F\left(\keps,\ue\right)_x
   + H\left(\keps,\leps,\ue\right)_y
   =  \Ieps_{1} + \Ieps_{2}+\Ieps_{3}+\Ieps_{4}+\Ieps_{5},
\end{align*}
where
\begin{align*}
	\Ieps_{1}&= \left(\eps \ue_x  f_u(\keps,\ue)\right)_x
	+  \left(\eps \ue_y  f_u(\keps,\ue)\right)_y,\\
	\Ieps_{2}&= - \eps \left(\ue_x\right)^2f_{uu}(\keps,\ue) 
	- \eps \left(\ue_y\right)^2f_{uu}(\keps,\ue), 
	\\ 
	\Ieps_{3} &= - \eps \ue_x f_{uk} (\keps,\ue) \keps_x 
	- \eps \ue_y f_{uk} (\keps,\ue) \keps_y, 
	\\ \Ieps_{4} &  = F_{k}(\keps,\ue)\keps_x
	+ H_{k}(\keps,\leps,\ue)\keps_y
	+ H_{l}(\keps,\leps,\ue)\leps_y
	\\ \Ieps_{5}
	&  = - f_{u}(\keps,\ue)f_k(\keps,\ue)\keps_x
	- f_{u}(\keps,\ue)g_l(\leps,\ue)\leps_y.
\end{align*}
Hence, there is a natural decomposition of $\Leto$ into six parts. We name 
the six parts $\Letonull$, $\Letoen$, $\Letoto$, $\Letotre$, 
$\Letofire$, and  $\Letofive$.

Regarding $\Letonull$,
$$
\abs{\Bigl\langle \Letonull,\phi\Bigr\rangle}
\le \tilde{C} \norm{\ue_{t}}_{L^\infty(\R_{+};L^1(\R^2))}
\norm{\phi}_{L^\infty(\Omega)}, \qquad \phi\in
C_0(\Omega),
$$
which yields $\norm{\Letonull}_{\CM(\Omega)}\le C$, where 
$\CM(\Omega)=\bigl(C_c(\Omega)\bigr)^\star$
denotes the space of bounded measures on $\Omega$.

In view of Lemma \ref{lem:energy} and the uniform boundedness 
of the solutions, it follows that 
$\left\{\Letoen\right\}_{\eps,\delta>0}$ 
is compact (and in fact converges to zero) 
in $W^{-1,2}(\Omega)$ 
and $\norm{\Letoto}_{\CM(\Omega)}\le C$, for each fixed $t>0$.

Next, for any $\phi\in C_{c}(\Omega)$, observe that
\begin{align*}
   &\abs{\iint_{\R^2} \Ieps_{3} \phi\dx\dy}
   \\ &  
   \le C \left\{\iint_{\R^2} \abs{\eps\keps_x} \abs{\keps_x}\right\}^{\frac12}
   \left\{\iint_{\R^2} \eps\left(\ue_x\right)^2\dx\dy\right\}^{\frac12}
   \\ &\qquad 
    + C \left\{\iint_{\R^2} \abs{\eps\keps_y} \abs{\keps_y}\right\}^{\frac12}
   \left\{\iint_{\R^2} \eps\left(\ue_y\right)^2\dx\dy\right\}^{\frac12}
\end{align*}
The point here is to have $\eps$ and $\delta$ in 
balance, so that we can ensure 
$\abs{\eps\keps_x},\,\abs{\eps\keps_y}\le C$. More precisely, we have 
$\abs{\eps\keps_x},\,\abs{\eps\keps_y}\le C \frac{\eps}{\delta}$, and 
by choosing $\eps, \delta$ according to \eqref{eq:balance} 
we achieve this balance. Consequently, 
$$
\abs{\iint_{\R^2} \Ieps_{3} \phi \dx\dy}
\le C \norm{\phi}_{L^\infty(\Omega)},
$$
and thus $\norm{\Letotre}_{\CM(\Omega)}\le C$, for each fixed $t>0$.

Finally, using the $BV$ regularity of the coefficients 
and the boundedness of the solutions, 
$$
\abs{\iint_{\R^2} \Ieps_{4} \phi \dx\dy}
\le C \norm{\phi}_{L^\infty(\Omega)},
$$
and thus  $\norm{\Letofire}_{\CM(\Omega)}\le C$, for each fixed $t>0$. 

Similarly, $\abs{\iint_{\R^2} \Ieps_{5} \phi \dx\dy}
\le C \norm{\phi}_{L^\infty(\Omega)}$, and thus $\norm{\Letofive}_{\CM(\Omega)}\le C$, for 
each fixed $t>0$. 

Summarizing, we have shown  that the sequence of distributions
$\left\{\Le\right\}_{\eps,\delta>0}$ satisfies the following two properties:   
\{i\} each distribution  is the sum of two terms --- one is compact in $W^{-1,2}(\Omega)$ 
and the other one is bounded in $\CM(\Omega)$. In addition, Lemma \ref{lem:max} implies that
\{ii\} $\left\{\Le\right\}_{\eps,\delta>0}$
belongs to a bounded subset of $W^{-1,\infty}(\Omega)$.
We now appeal to Murat lemma \cite{Murat:Hneg}, which guarantees that 
$\left\{\Le\right\}_{\eps,\delta>0}$ belongs to a compact 
subset of $W^{-1,2}(\Omega)$. This concludes the proof of the first part of the 
lemma, since $\Omega$ was an arbitrary bounded open subset of $\R^2$. 
The second part of the lemma can be proved in a similar way.
\end{proof}

\begin{proof}[Concluding the proof of Theorem \ref{thm:main}]
By combining Lemmas \ref{lem:CC} and \ref{lem:Hneg_viscous}, we 
conclude that $\ue(\cdot,\cdot,t)$ is pre-compact a.e.~for each $t\in[0,T]$. Together 
with a diagonal argument, we can prove that $\ue(\cdot,\cdot,t)$ converges 
along a subsequence a.e.~in $\R^2$ and in 
$L^1_{\loc}(\R^2)$, for each fixed $t>0$. Lemma \ref{lem:time} implies that
$$
\norm{\ue(\cdot,\cdot,t+\tau)-\ue(\cdot,\cdot,t)}_{L^1(\R^2)}
\le C \tau, \qquad \forall \tau\in (0,T-\tau),
$$
and using this $L^1$ time continuity estimate it takes 
a standard density argument to show that there exists a subsequence of 
$\seq{\ue}_{\eps,\delta>0}$ that converges 
to a limit function $u$ a.e.~in $\R^2\times\R_{+}$ and in 
$L^1_{\loc}(\R^2\times\R_{+})$. Moreover, 
the limit $u$ belongs to $L^\infty(\R^2\times\R_{+})\cap 
\mathrm{Lip}(\R_{+};L^1(\R))$. 

Equipped with the strong convergence it is easy 
to prove that the limit $u$ is a weak solution.  
This concludes the proof of Theorem \ref{thm:main}. 
\end{proof}

We close the paper by a couple of remarks.
First, we note that the apriori bounds in Lemmas 
\ref{lem:max},  \ref{lem:time} and \ref{lem:energy}, being  uniform in time, enabled 
us to deduce pre-compactness at each fix $t>0$, thus circumventing the 
temporal argument required in \cite[Appendix A]{TRB}.
Second, we have herein exclusively dealt with problems that 
are spatially two-dimensional. A possible strategy for 
going beyond two dimensions is to apply the compactness 
framework of Panov \cite{Panov:1994gf,Panov:1999mz,Panov:1995ys}.

\def\cprime{$'$} \def\cprime{$'$} \def\cprime{$'$}

\end{document}